\documentclass[12pt, letterpaper]{article}
\usepackage{amsfonts}
\usepackage{amssymb}
\usepackage{latexsym}  

\input cyracc.def
\newfam\cyrfam

\font\twelvecyr=wncyr10 scaled 1200 
\font\twelveitcyr=wncyi10 scaled 1200 
\font\twelvebfcyr=wncyb10 scaled 1200 
\def\cyr{\fam\cyrfam\twelvecyr\cyracc}
\def\itcyr{\fam\cyrfam\twelveitcyr\cyracc}
\def\bfcyr{\fam\cyrfam\twelvebfcyr\cyracc}

\newcommand\ZZZ{\mathbb{Z}}

\newcommand\rest{\upharpoonright}

\newcommand\Atop{\mathrm{Atop}}
\newcommand\Sym{\mathrm{Sym}}
\newcommand\Aut{\mathrm{Aut}}
\newcommand\NAut{\mathit{N}\!\mathrm{Aut}}
\newcommand\fix{\mathrm{fix}}

\newcommand\LL{\mathcal{L}}

\newcommand\TT{\mathcal{T}}

\newcommand\UU{\mathcal{U}}
\newcommand\VV{\mathcal{V}}

\newcommand\RCC{\mathit{RCC}}
\newcommand\LCC{\mathit{LCC}}

\newcommand\iv{^{-1}} 

\newcommand\eop{ ${\vcenter
   {\hrule
   \hbox{\vrule height 7pt \kern 7pt \vrule height 7pt}
   \hrule}}$}

\newcommand\ld{\backslash}
\newcommand\rd{/}

\newcommand\Hom{\mathrm{Hom}}

\newcommand\sbl[1]{\langle#1\rangle}   

\newtheorem{theorem}{Theorem}[section]
\newtheorem{corollary}[theorem]{Corollary}
\newtheorem{definition}[theorem]{Definition}
\newtheorem{lemma}[theorem]{Lemma}
\newtheorem{proposition}[theorem]{Proposition}
\newtheorem{remark}[theorem]{Remark}

\newenvironment{itemizz} %
{\begin{itemize}\setlength{\itemsep}{-1mm}}{\end{itemize}}

\newenvironment{proof}{\noindent\textit{Proof.}}{\eop\medskip}

\begin{document}

\title{Diassociativity in Conjugacy Closed Loops}

\author{Michael K. Kinyon,
Kenneth Kunen\thanks
{Author supported by NSF Grant DMS-0097881},
and J.D. Phillips}

\maketitle

\begin{abstract}
Let $Q$ be a conjugacy closed loop, and $N(Q)$ its nucleus.
Then $Z(N(Q))$ contains all associators of elements of $Q$.
If in addition $Q$ is diassociative (i.e., an extra loop),
then all these associators have order $2$.
If $Q$ is power-associative and $|Q|$ is finite and
relatively prime to $6$, then $Q$ is a group.
If $Q$ is a finite non-associative extra loop,
then $16 \mid |Q|$.
\end{abstract}


\section{Introduction}
\label{sec-introduction}

The notion of a conjugacy closed loop (CC-loop) is due to
Goodaire and Robinson \cite{GRA},
and independently to {\cyr So\u\i kis} \cite{SO},
with somewhat different terminology.  Following, approximately, \cite{GRA}:

\begin{definition}
\label{def-cc-alt}
A loop $(Q, \cdot)$ is {\em conjugacy closed} (or a {\em CC-loop}) if
and only if there are functions $f,g :  Q\times Q \to Q$ such that for all $x,y,z$:
\[
\RCC:\  x \cdot yz = f(x,y)\cdot x z  \qquad\qquad
\LCC:\  zy \cdot x = zx \cdot g(x,y)   \ \ .
\]
\end{definition}

As usual, define the left and right multiplications by
$xy = x R_y = y L_x$, so that $R_y$ and $L_x$ are permutations of
the set $Q$.  Using these, we can express ``CC-loop'' in terms
of conjugations:

\begin{lemma}
\label{lemma-conj}
A loop $Q$ is a CC-loop if and only if there exist functions
$f,g: Q\times Q\to Q$ such that
\[
L_x\iv L_y L_x =L_{f(x,y)}
\qquad and \qquad
R_x\iv R_y R_x =R_{g(x,y)}.
\]
\end{lemma}
\begin{proof}
$\RCC$ and $\LCC$ assert that
$L_y L_x =L_x L_{f(x,y)}$ and $R_y R_x =R_x R_{g(x,y)}$.
\end{proof}

Thus, in a CC-loop, the left multiplications are closed under
conjugation and the right multiplications are closed under
conjugation; hence the name ``conjugacy-closed''.

These loops have a number of interesting properties,
surveyed in Sections \ref{sec-back} and \ref{sec-nuc};
for example, by \cite{GRA}, the left and right inner mappings are
automorphisms.  These properties allow a rather detailed structural
analysis to be made; in particular, all CC-loops of orders
$p^2$ and  $2p$ (for primes $p$) are known (see \cite{KUNC}).
This paper yields additional structural information about CC-loops ---
especially for the ones which are
\textit{power-associative} (that is, each $\sbl{x}$ is a group) or
\textit{diassociative} (that is, each $\sbl{x,y}$ is a group).

It is shown in \cite{GRB} that the CC-loops which are diassociative
(equivalently, Moufang) are the extra loops studied by Fenyves
\cite{FENA, FENB}. By \cite{FENB}, if $Q$ is an extra loop, then
$Q/N(Q)$ is a boolean group (where $N(Q)$ is the nucleus).  It is
immediate that a finite extra loop of odd order is a group.  We show
here (Corollary \ref{cor-finite-sixth})
that a finite power-associative CC-loop of order relatively prime
to $6$ is a group. The ``6'' cannot be improved, since there
are non-associative power-associative CC-loops of order $16$
(e.g., the Cayley loop) and of order $27$ (see Section \ref{sec-ex})
(we do not know if there are ones of
order divisible by $6$ but not by $4$ or $9$).   Also, one cannot drop
the ``power-associative'', since by \cite{GRA},
there are non--power-associative CC-loops of order
$p^2$ for every odd prime $p$.

More generally, we show that every power-associative CC-loop satisfies
a weakening of diassociativity --- namely, $\sbl{x,y}$ is a group
whenever $x$ is a cube and $y$ is a square.  Then, if $|Q|$
is relatively prime to $6$, every element must be a sixth power
by the Lagrange property, so that $Q$ is diassociative, and
hence an extra loop of odd order, and hence a group.  Of course,
we must verify that the Lagrange property really holds for CC-loops,
since it can fail for loops in general.  This is easy to do
(see Corollary \ref{cor-lagrange}) using
the result of {\cyr Basarab} \cite{BAS}.  He showed that for any CC-loop,
$Q/N(Q)$ is an abelian group (this answers a question
from \cite{GRA}); we include a proof of this here (see Theorem
\ref{thm-bas}), since it is fairly short using the notion of autotopy
(see {\cyr  Belousov} \cite{BEL} II\S3 or Bruck \cite{BR2} VII\S2),
together with some facts about the autotopies of CC-loops proved
by Goodaire and  Robinson (see \cite{GRA} and Section \ref{sec-back}).

We also establish two theorems about general CC-loops.
First, whenever $S \subseteq Q$ and $S$ \textit{associates} in the sense
that $x \cdot yz = xy \cdot z$ holds for all $x,y,z\in S$,
we prove that $\sbl{S}$ also associates, and hence is a group
(see Corollary \ref{cor-assoc}).
Second (see Theorems  \ref{thm-group2} and \ref{thm-group1}),
we use this fact to show that 
$\sbl{b, c^2}$ and $\sbl{b^2, c}$ are groups whenever
$\sbl{b}$ and $\sbl{c}$ are groups and $c$ satisfies
$ c \cdot ( (xc)\ld 1) = x \ld 1$  (such $c$ are called
\textit{WIP elements}; see Definition~\ref{def-WIP}).

Finally, in a power-associative CC-loop,
we show that all cubes are WIP elements (see Section \ref{sec-pow}),
so that the subloop generated by a square and a cube always is a group.

Our investigations were aided by the computer programs
OTTER, developed by McCune \cite{MC}, and SEM,
developed by J. Zhang and H. Zhang \cite{ZZ}. 


\section{Background}
\label{sec-back}

Let $Q$ be a loop. We shall reformulate the notion of CC-loop in terms
of autotopisms, the definition of which we now recall.

\begin{definition}
\label{def-Atop}
Let $\Sym(Q)$ denote the group of all permutations of the set $Q$,  and
let $I$ denote the identity element of $\Sym(Q)$. A triple
$(\alpha,\beta,\gamma)\in (\Sym(Q))^3$ is an \emph{autotopism} of
$Q$ if $y\alpha \cdot z\beta = (yz)\gamma$ for all $y,z \in Q$.
Let $\Atop(Q)$ denote the set of all autotopisms of $Q$.
\end{definition}

It is easy to see that $\Atop(Q)$ is a subgroup of $(\Sym(Q))^3$.

\begin{lemma}
\label{lemma-CCAtop}
A loop $Q$ is a CC-loop if and only if there exist mappings
$F,G: Q\to \Sym(Q)$ such that
\[
 (F_x , L_x, L_x)
\qquad \mathrm{and}\qquad
 (R_x, G_x, R_x)
\]
are in $\Atop(Q)$.  In this case, $F_x$ and $G_x$ are given by:
$yF_x = f(x,y)$ and $yG_x = g(x,y)$ (see Definition \ref{def-cc-alt}).
\end{lemma}

We shall also use the division and the left and right inverse permutations:

\begin{definition}
In any loop $Q$, define permutations $\rho$ and $\lambda$, along
with $D_x$ for $x\in Q$, by:
$$
y \lambda = 1 \rd y \qquad
y \rho = y \ld 1 \qquad
y D_x = y \ld x \ \ .
$$
We write $y^\lambda, y^\rho$ for $y \lambda, y \rho$, respectively; when these
values are the same, they are denoted by $y\iv$.
If $y^\lambda =  y^\rho$ for \emph{all} $y$, we let $J = \lambda = \rho$.
\end{definition}

Note that $y D_x\iv = x \rd y$,  $\rho = D_1$, and $\lambda = D_1\iv$.
These permutations are used in the following explicit expressions
for $F_x$ and $G_x$, which are obtained from Definition \ref{def-cc-alt}:

\begin{lemma}
\label{lemma-getfg}
For all $z$,
\[
F_x = R_z L_x R_{xz}\iv = D_z L_x D_{xz}\iv
\qquad and \qquad
G_x = L_z R_x L_{zx}\iv = D_z\iv R_x D_{zx} .
\]
In particular,
\[
\begin{array}{rcccccccl}
f(x,y) &=& (xy) \rd x &=& x \cdot yx^{\rho} &=& x \rd (xy^{\rho}) &=& [x(y\ld x^{\rho})]^{\lambda}\\
F_x &=& L_x R_x\iv &=& R_{x^{\rho}} L_x &=& \rho L_x D_x\iv &=& D_{x^{\rho}} L_x\, \lambda \\
\\
g(x,y) &=& x \ld (yx)  &=& x^{\lambda}y \cdot x &=& (y^{\lambda}x) \ld x &=& [(x^{\lambda}\rd y)x]^{\rho} \\
G_x &=& R_x L_x\iv &=& L_{x^{\lambda}} R_x &=& \lambda R_x D_x &=& D_{x^{\lambda}}\iv R_x \,\rho
\end{array}
\]
\end{lemma}

\begin{proof}
$F_x = R_z L_x R_{xz}\iv$ is immediate from $\RCC$.
Replacing the $z$ in $\RCC$ by $y\ld z$ we obtain
$xz = f(x,y)\cdot x(y\ld z)$, which yields $F_x =  D_z L_x D_{xz}\iv$.
The rest of the expressions for $F_x$ are obtained by setting $z$
to equal either $1$ or $x^{\rho}$.  The expressions for $G_x$ are
likewise obtained from $\LCC$.
\end{proof}

\begin{corollary}
\label{coro-fginverses}
$F_x G_x = G_x F_x = I$
\end{corollary}

\begin{corollary}
\label{coro-preAAIP}
$(x \cdot yx^{\rho}) \cdot xy^{\rho} = x$
and
$y^{\lambda}x \cdot (x^{\lambda}y\cdot x) = x$.
\end{corollary}

We remark that a loop has the \textit{anti-automorphic inverse property (AAIP)}
iff it satisfies $(xy)^{\rho} = y^{\rho}x^{\rho}$. This is equivalent to
$(xy)^{\lambda} = y^{\lambda}x^{\lambda}$, and these conditions imply
$\rho = \lambda$. In Corollary \ref{coro-preAAIP}, $yx^{\rho}$ can be an
arbitrary element of the loop, so the AAIP would give us $xz\cdot z\iv = x$,
which is the inverse property (IP).
Since a CC-loop with the IP is an extra loop \cite{GRB}, we have:

\begin{remark}
\label{remark-aaip}
A CC-loop with the AAIP is an extra loop.
\end{remark}

The following lemma lists some
additional conjugation relations among the
left and right translations; (3) and (4) are from
\cite{KUNC}, Lemma 3.1:

\begin{lemma}
\label{lemma-LR-conj}
In any CC-loop:
\[
\begin{array}{rlcl}
1. & L_x L_y L_x\iv = L_{g(x,y)}  & &
    R_x R_y  R_x\iv =  R_{f(x,y)} \\
2. & x \cdot g(x,y) z = y \cdot x z & &
   z f(x,y) \cdot x = zx \cdot y \\
3. & L_x\iv R_y L_x = R_x\iv R_{xy} = R_{y\rd x^{\rho}} R_x\iv & &
   R_x\iv L_y R_x = L_x\iv L_{yx} = L_{x^{\lambda} \ld y} L_x\iv \\
4. & L_x R_y L_x\iv = R_{x^{\rho}}\iv R_{x \ld y} = 
            R_{yx^\rho} R_{x^{\rho}}\iv & &
   R_x L_y R_x\iv = L_{x^{\lambda}}\iv L_{y \rd x} = 
            L_{x^\lambda y} L_{x^{\lambda}}\iv 
\end{array}
\]
\end{lemma}

\begin{proof}
For (1), use Lemma \ref{lemma-conj} and Corollary \ref{coro-fginverses}.
(2) is equivalent to (1).
For the first equality of (3), use Lemma \ref{lemma-getfg}
and Corollary \ref{coro-fginverses} to get
$$
R_x L_x\iv R_z L_x = G_x R_z L_x = F_x\iv R_z L_x = R_{xz} \ \ .
$$
For the second one, use (1) and Lemma \ref{lemma-getfg} to get
$ R_x R_{y\rd x^{\rho}} R_x\iv = R_{f(x,y\rd x^{\rho})} = R_{xy}$.
For the first equality of (4), use Lemma \ref{lemma-getfg} with $z$ replaced
by $x\ld z$ to obtain $R_{x^{\rho}} L_x R_z = R_{x\ld z} L_x$. For
the second equality, use Lemmas \ref{lemma-conj} and \ref{lemma-getfg}:
$R_{x^{\rho}}\iv  R_{x\ld y} R_{x^{\rho}} =
R_{g(x^{\rho}, x\ld y)} = R_{yx^\rho}$.
\end{proof}

The left nucleus ($N_\lambda$), the middle nucleus ($N_\mu$), the
right nucleus ($N_\rho$), and the nucleus ($N$) are defined by:

\begin{definition}
\label{def-nuc}
Let $Q$ be a loop.
\[
\begin{array}{rcl}
N_\lambda(Q) &:=& \{ a\in Q : 
  \forall x, y \in Q \,[\, a\cdot xy = ax\cdot y \,] \} \\
N_\mu(Q) &:=& \{ a\in Q : 
  \forall x, y \in Q \,[\, xa\cdot y = x\cdot ay \,] \} \\
N_\rho(Q) &:=& \{ a\in Q : 
  \forall x, y \in Q \,[\, x\cdot ya = xy\cdot a \,] \} \\
N(Q) &:=& N_\lambda(Q) \cap  N_\mu(Q) \cap N_\rho(Q)
\end{array}
\]
\end{definition}

It is easy to verify the following equivalents, in terms of autotopy.

\begin{lemma}
\label{lemma-nuc-atop}
For any loop $Q$:
\begin{itemizz}
\item [1.]
$N_\lambda(Q) = \{a \in Q : (L_a, I, L_a) \in \Atop(Q)\}$. \\
$N_\mu(Q) = \{a \in Q : (R_a\iv, L_a, I) \in \Atop(Q)\}$. \\
$N_\rho(Q) = \{a \in Q : (I, R_a, R_a) \in \Atop(Q)\}$.
\item [2.]
If $(\alpha,I, \gamma) \in \Atop(Q)$, then $\alpha = \gamma$,
$1\alpha \in N_\lambda(Q)$ and $\alpha = L_{1\alpha}$.
\end{itemizz}
\end{lemma}
\begin{proof}
For (2), $x\alpha \cdot y = (xy)\gamma$, so taking $y=1$ gives
$\alpha = \gamma$. Then let $a = 1\alpha$ and $x=1$ to
obtain $a y = y\alpha$, so that $ax \cdot y = a\cdot xy$.
\end{proof}

\begin{definition}
For any loop $Q$,
$Z(Q) = \{x \in N(Q) : \forall y (xy = yx)\} $,
\end{definition}

By Goodaire and Robinson \cite{GRA}:

\begin{theorem}
\label{thm-nuc-basic}
In any CC-loop $Q$, $N(Q) =  N_\lambda(Q) =  N_\mu(Q)  = N_\rho(Q)$
and $N(Q)$ is a normal subloop of $Q$.
Also, $Z(Q) =\{x \in Q : \forall y (xy = yx)\} $, so that
every commutative CC-loop is a group.
\end{theorem}

Autotopies are useful for producing automorphisms:

\begin{lemma}
\label{lemma:Aut-Atop}
In any loop $Q$, if $1 \alpha = 1$ and either
$(\alpha, \beta, \alpha) \in \Atop(Q)$ or
$(\beta, \alpha, \alpha) \in \Atop(Q)$, then $\alpha = \beta$
and $\alpha$ is an automorphism.
\end{lemma}

\begin{proof}
If $(\alpha, \beta, \alpha) \in \Atop(Q)$,
we have $x\alpha \cdot y\beta = (xy)\alpha$.
Setting $x = 1$ yields $\alpha = \beta$.
\end{proof}

Following Bruck \cite{BR2} \S IV.1, define the generators of the
right and left inner mapping groups by:

\begin{definition}
$R(x,y) := R_x R_y R_{xy}\iv$ and $L(x,y) := L_x L_y L_{yx}\iv$.
\end{definition}

Using Lemmas \ref{lemma-CCAtop} and \ref{lemma:Aut-Atop}, we get

\begin{lemma}[\cite{GRA}]
\label{lemma-aut}
In any CC-loop, $R(x,y)$ and $L(x,y)$ are automorphisms
for all $x,y$.
\end{lemma}

The following definitions will be useful in
Sections \ref{sec-wip} and \ref{sec-pow}:

\begin{definition}
\label{defn-E}
$E_x = R(x, x^{\rho}) = R_x R_{x^\rho} $.
\end{definition}

\begin{definition}
\label{def-PA}
$a$ is a \emph{power-associative} element if $\sbl{a}$ is a group. A loop 
is \emph{power-associative} iff every element is power-associative.
\end{definition}

By (\cite{KUNC}, Lemma 3.20):

\begin{lemma}
\label{lem-PA}
Let $Q$ be a CC-loop. The following are equivalent for an element $a\in Q$:
(i) $a$ is power-associative; (ii) $1 \rd a = a \ld 1$; (iii) $a\cdot aa = aa\cdot a$.
In particular, $Q$ is power-associative if and only if $\rho = \lambda$.
\end{lemma}

Following Osborn \cite{OS}:

\begin{definition}
\label{def-WIP}
$c$ is a \emph{WIP element} (briefly: $c$ is WIP) iff $\lambda R_c  \rho = L_c\iv$.
A loop has the \emph{weak inverse property} iff every element is WIP;
\end{definition}

For convenience, we collect the following easy equivalents of WIP.

\begin{lemma}
\label{lemma-wip-equiv}
In any loop, each of the following four equations is equivalent
to the statement that $c$ is a WIP element:
$$
\begin{array}{ccc}
\lambda R_c  \rho = L_c\iv &\qquad & \rho L_c \lambda = R_c\iv  \\
R_c \rho L_c = \rho  &\qquad & L_c \lambda R_c = \lambda
\end{array}
$$
\end{lemma}


\section{$Q/N(Q)$}
\label{sec-nuc}
\begin{theorem}[{\bfcyr  Basarab} \cite{BAS}]
\label{thm-bas}
For a CC-loop $Q$,  $Q/N(Q)$ is an abelian group.
\end{theorem}
\begin{proof}
By Goodaire and Robinson \cite{GRA},
every CC-loop is a G-loop;
that is, it is isomorphic to all its loop isotopes.
In particular, for any element $v$,
the isotope $(Q; \circ)$ defined
by $ x \circ y = x \cdot (v \ld y)$ satisfies $\RCC$:
$$
x \cdot (v \ld (y \cdot (v \ld z))) =
h(x,y,v)\cdot (v \ld (x \cdot (v \ld z)))  \ \ ,
$$
where $h : Q^3 \to Q$.  Replacing $z$ by $vz$, this simplifies to:
$$
x \cdot (v \ld (y \cdot  z)) =
h(x,y,v)\cdot (v \ld (x \cdot z))  \ \ .
$$
We may set $z = 1$ to get
$h(x,y,v) = (x  (v \ld y) ) \rd (v \ld x)$, so we have
$$
x \cdot (v \ld (y \cdot  z)) =
[(x  (v \ld y) ) \rd (v \ld x) ] \cdot [v \ld (x \cdot z)]  \ \ ,
$$
which implies that 
$ ( L_v\iv L_x R_{v\ld x}\iv ,\; L_x L_v\iv ,\; L_v\iv L_x )  
\in \Atop(Q)$ for all $x$ and $v$.
Since also $( F_v F_x\iv,\; L_v L_x\iv ,\; L_v L_x\iv ) \in \Atop(Q)$
by Lemma \ref{lemma-CCAtop}, we have
$$
( F_v F_x\iv L_v\iv L_x R_{v\ld x}\iv ,\; I,\;  L_v L_x\iv L_v\iv L_x )
 \in \Atop(Q)
$$
Then, by Lemma \ref{lemma-nuc-atop},
$1\, F_v F_x\iv L_v\iv L_x R_{v\ld x}\iv = (x (v \ld 1))\rd (v \ld x  ) 
\in N_\lambda(Q)$.
Applying Theorem \ref{thm-nuc-basic}, 
$Q / N(Q)$ is a CC-loop satisfying the additional equation
$(x (v \ld 1))\rd (v \ld x  ) = 1$, or $x v^\rho  = v \ld x$.
Hence, in $Q/N(Q)$, we have (using Lemma \ref{lemma-getfg})
$f(v,y) = v \cdot y v^\rho = y$, so that $\RCC$ becomes
$x \cdot yz = y\cdot x z$.
Setting $z = 1$, we get $xy = yx$, so that $Q/N(Q)$ is commutative
and satisfies the associative law,
$x \cdot zy = x z \cdot y $.
\end{proof}

This is roughly the proof in \cite{BAS}, although {\cyr Basarab} studies
in more detail those loops $Q$ such that
$Q$ and all its loop isotopes satisfy $\RCC$.

Recall that a finite loop has the \textit{weak Lagrange property} if the order
of any subloop divides the order of the loop, and a finite loop has the
\textit{strong Lagrange property} if every subloop has the weak Lagrange
property \cite{PF}.  In general if $H$ is a normal subloop of $Q$,
and $H$ and $Q/H$ both have the strong Lagrange property, then so does $Q$
(see Bruck \cite{BR2}, \S V.2,  Lemma 2.1; see also \cite{CKRV}).
It is now immediate from Theorem \ref{thm-bas} that:

\begin{corollary}
\label{cor-lagrange}
Every finite CC-loop has the strong Lagrange property.
\end{corollary}

\begin{corollary}
\label{cor-orders}
If $Q$ is a finite power-associative CC-loop and
$|Q|$ is relatively prime to $n$, then every element of
$Q$ is an $n^{\mathrm{th}}$ power.
\end{corollary}

The following Cauchy property is also immediate from Theorem \ref{thm-bas}:

\begin{corollary}
If $Q$ is a finite power-associative CC-loop and
$|Q|$ is divisible by a prime $p$, then $Q$ contains
an element of order $p$.
\end{corollary}

Also, the fact that finite $p$-groups have non-trivial centers generalizes
to:

\begin{corollary}
\label{cor-center}
If $Q$ is a CC-loop of order $p^n$ for some prime $p$ and $n > 0$,
then
\begin{itemizz}
\item[1.] $|Z(Q)| = p^r$, where $ r \ne 0$ and $r \ne n-1$.
\item[2.] For all $m \le n$, $Q$ has a normal subloop of order $p^m$.
\end{itemizz}
\end{corollary}
\begin{proof}
For (1):
Let $N$ be the nucleus.  For $x \in Q$, let $T_x = R_x L_x\iv$.
By \cite{GRA}, each $T_x\rest N$ is an automorphism of $N$.
Furthermore, if we define $\TT : Q \to \Aut(N)$ by 
$\TT(x) = T_x\rest N$, then $\TT$ is a homomorphism
by \cite{KUNC}, Corollary 3.7.
Thus, $\TT(Q)$ is a subgroup of $\Sym(N)$, 
and $|\TT(Q)|$ is a power of $p$, so the size of each orbit is
a power of $p$.  Since $|N| = p^\ell$ for some $\ell > 0$,
there must be at least $p$ elements $y$ whose orbit is a singleton
(equivalently, $y \in Z(Q)$).
Hence $r \ne 0$.

If $r \ge n-1$, then $Q = \sbl{Z(Q) \cup \{a\}}$ for any
$a \ne Z(Q)$, but then $Q$ is commutative, so $r = n$.

For (2):  Let $P$ be a subgroup of $Z(Q)$ of order $p$.
Then the $m = 1$ case is immediate, using $P$, and the
case $1 < m \le n$ follows by applying induction to $Q/P$.
\end{proof}


\section{Associators and Inner Mappings}
\label{sec-associators}

\begin{definition}
\label{defn-associator}
In a loop $Q$, \emph{associators} are denoted by:
$$
(x,y,z) := (x \cdot yz)\, \ld\, (xy \cdot z)  \qquad
[x,y,z] := (x \cdot yz)\, \rd\, (xy \cdot z) \ \ .
$$
\end{definition}

Since the two notions of ``associator'' are mirrors of each other,
we concentrate on $(x,y,z)$ in the following:

\begin{lemma}
\label{lemma-nuc-assoc}
In any loop $Q$ with nucleus $N$, if $a\in N$, then
\begin{itemizz}
\item[(i)] $(ax,y,z) = (x,y,z)$
\item[(ii)] $(xa,y,z) = (x,ay,z)$
\item[(iii)] $(x,ya,z) = (x,y,az)$
\item[(iv)] $(x,y,za) = a\iv (x,y,z)a$
\end{itemizz}
In addition, if $N$ is normal in $Q$, then
\begin{itemizz}
\item[(v)] $(xa,y,z) = (x,y,z)$
\item[(vi)] $(x,ya,z) = (x,y,z)$
\item[(vii)] $a\iv (x,y,z) a = (x,y,z)$
\end{itemizz}
\end{lemma}

\begin{proof}
(i)-(iv) are straightforward consequences of the definitions.
Now assume $N$ is normal in $Q$. Then for $u\in Q$,
$ua = bu$ for some $b \in N$. Thus (v) follows from (i),
(vi) follows from (ii) and (v), and (vii) follow from (iv),
(iii), and (vi).
\end{proof}

Theorem \ref{thm-bas} implies that associators are nuclear,
so we have:

\begin{corollary}
\label{coro-nuc-center}
The nucleus of a nonassociative CC-loop
has a nontrivial center which contains the subgroup
generated by the associators.
\end{corollary}

\begin{theorem}
\label{thm-perm}
In a CC-loop, the associators $(x,y,z)$ and $[x,y,z]$
are invariant under all permutations of the set $\{x,y,z\}$.
\end{theorem}

\begin{proof}
It is enough to prove that $(x,y,z) = (y,x,z)$ and
$(x,y,z) = (x,z,y)$, since the transpositions
$(x \; y)$ and $(y \; z)$ generate $\Sym(\{x,y,z\})$.

For $(x,y,z) = (y,x,z)$:
$x \cdot yz = f(x,y)\cdot x z$ by $\RCC$
and
$xy \cdot z = f(x,y) x \cdot z$ by 
Lemma \ref{lemma-getfg}, so  $(x,y,z) = (f(x,y),x,z)$.
By Theorem \ref{thm-bas}, there exists $a \in N$ such that $f(x,y) = ay$,
so $(f(x,y),x,z) = (ay, x,z) = (y, x, z)$ by Lemma \ref{lemma-nuc-assoc}(i).

For $(x,y,z) = (x,z,y)$:  Apply a similar argument, using $\LCC$.
\end{proof}

\begin{lemma}
\label{lemma-LR-assoc}
In a CC-loop,
$$
zL(y,x) = z (x,y,z)\iv  \qquad
xR(y,z) =  [x,y,z]\iv x  \ \ .
$$
\end{lemma}
\begin{proof} We have:
$$
(x \cdot yz) (x,y,z) = (xy \cdot z)  \qquad
[x,y,z] (xy \cdot z) = (x \cdot yz)  \ \ .
$$
Since associators are nuclear, this can be rewritten as
$$
\{z(x,y,z)\} L(y,x) = z \qquad
\{[x,y,z] x\} R(y,z) = x \ \ .
$$
Now use the fact that $L(y,x)$ and $R(y,z)$ are automorphisms
which fix all elements of the nucleus.
\end{proof}

Applying Theorem \ref{thm-perm}:

\begin{corollary}
\label{coro-inner-perm}
In a CC-loop, $L(x,y) = L(y,x)$ and $R(x,y) = R(y,x)$.
\end{corollary}

Furthermore, the $L(x,y)$ and $R(u,v)$ all commute with 
each other; more generally, they commute with all nuclear automorphisms:

\begin{definition}
\label{defn-nuclear}
Let $Q$ be a loop with nucleus $N$. An automorphism
$\alpha$ of $Q$ is \emph{nuclear} iff
$x\alpha \in xN$ for each $x\in Q$.
$\; \NAut(Q)$ is the set of nuclear automorphisms of $Q$.
\end{definition}

\begin{lemma}
$\NAut(Q)$ is a normal subgroup of $\Aut(Q)$.
\end{lemma}

\begin{theorem}
\label{thm-commutes}
Let $Q$ be a CC-loop.  Then $Z(\NAut(Q))$ contains
all $R(x,y)$ and $L(x,y)$.
\end{theorem}
\begin{proof}
The $R(x,y)$ and $L(x,y)$ are automorphisms by 
Lemma \ref{lemma-aut} and nuclear by Theorem \ref{thm-bas}.
Now, if $\alpha$ is nuclear, we have
$$
 z L(y,x) \alpha = \{z (x,y,z) \iv \}\alpha =
z \alpha \cdot (x a, y b, z c) \iv = z \alpha \cdot (x,y,z) \iv \ \ ,
$$
where $a,b,c \in N(Q)$, whereas
$$
   z \alpha L(y,x) = ( z \alpha ) \cdot  (x, y, z \alpha)\iv
     = ( z \alpha )  \cdot (x, y, z d)\iv = z \alpha  \cdot (x,y,z)\iv \ \ ,
$$
where $d\in N(Q)$.
\end{proof}

\begin{corollary}
\label{coro-commute}
In a CC-loop, the group generated by all the automorphisms
$R(x,y)$ and $L(x,y)$ is abelian.
\end{corollary}

We conclude this section with some applications to extra loops.
As mentioned in the Introduction, extra loops are Moufang
CC-loops, and as CC-loops, they have several characterizations.
Indeed, each of the following is sufficient for a CC-loop
to be an extra loop: (i) the left or right alternative laws
($x\cdot xy = x^2y$ or $xy\cdot y = xy^2$),
(ii) the flexible law ($x\cdot yx = xy\cdot x$) \cite{GRB},
(iii) the AAIP (see Remark \ref{remark-aaip}),
(iv) the left or right IP ($x\ld y = x^{\lambda}y$
or $x\rd y = xy^{\rho}$) \cite{GRB}, (v) diassociativity.
The nucleus of an extra loop contains every square \cite{FENB}.
However, there are non-extra CC-loops $Q$ in which all squares
are in the nucleus; $Q$ can both be power-associative and have
the weak inverse property; see Section \ref{sec-ex}.

\begin{lemma}
\label{lemma-assoc-inv}
In a CC-loop, $z (x,y,z)\iv = (x,y,z^{\lambda}) z$.
\end{lemma}

\begin{proof}
Applying the automorphism $L(y,x)$ to the equation $z^\lambda z = 1$,
and using Lemma \ref{lemma-LR-assoc}, we get
$ z^{\lambda} (x,y,z^{\lambda})\iv \cdot z (x,y,z)\iv =  1 = z^\lambda z$.
The result now follows because associators are in the nucleus.
\end{proof}

\begin{lemma}
\label{lemma-square-nuc}
Let $Q$ be a CC-loop such that $N(Q)$ contains every square.
For $i = 1,2,3$, choose $\epsilon_i \in \{I, \rho, \lambda\}$.
Then $(x,y,z) = (x\epsilon_1 ,y\epsilon_2 ,z\epsilon_3)$.
Hence $L(y,x) = L(y\epsilon_1 , x\epsilon_2)$.
\end{lemma}

\begin{proof}
Note that $z^2 z^{\lambda} = z^{\rho} z^2 = z$
(since $z^2 z^{\lambda} z = z^2$).  Then,
Lemma \ref{lemma-nuc-assoc} implies
$(x,y,z) = (x,y,z^{\lambda}) = (x,y,z^{\rho})$. The
remainder follows from Theorem \ref{thm-perm} and
Lemma \ref{lemma-LR-assoc}.
\end{proof}

\begin{theorem}
\label{thm-extra}
Let $Q$ be an extra loop.
\begin{itemizz}
\item[1.] $L(x,y) = R(x,y) = L(y,x) = R(y,x)$ and $L(x,y)^2 = I$.
\item[2.] $(x,y,z) = [x,y,z]$  and $(x,y,z)^2 = 1$,
so that the subgroup of $N(Q)$ generated by the associators
is a boolean group.
\item[3.] Each $(x,y,z)$ commutes with $x$, $y$, and $z$.
\end{itemizz}
\end{theorem}

\begin{proof}
(1) In Moufang loops, $R(x\iv, y\iv) = L(x,y) = L(y,x)\iv$
(see \cite{BR2}, Lemma VII.5.4).  Now apply
Corollary \ref{coro-inner-perm} and Lemma \ref{lemma-square-nuc}.

(2) In diassociative loops, $(x,y,z)\iv = [z\iv, y\iv,x\iv]$.
Now, apply (1), along with  Theorem \ref{thm-perm} and 
Lemmas \ref{lemma-LR-assoc} and \ref{lemma-square-nuc}.

(3) This follows from (2), Lemmas \ref{lemma-assoc-inv} and
\ref{lemma-square-nuc}, and Theorem \ref{thm-perm}.
\end{proof}

Hence, the nucleus of a nonassociative extra loop 
must contain elements of order 2.

\begin{corollary}
If $Q$ is a finite nonassociative extra loop, then $16 \mid |Q|$.
\end{corollary}
\begin{proof}
Since the order of $N = N(Q)$  is even, and $Q/N$ is a boolean group,
it is sufficient to show that $|Q:N| \ge 8$.
Choose $a \notin N = N_\mu$, and then choose $b$ such that
$R(a,b) \ne I$ (that is, $(xa)b \ne x(ab)$ for some $x$).
Then $N < \fix(R(a,b)) < Q$, since $a,b \in \fix(R(a,b))$.
Next, note that $\sbl{N\cup\{a\}} = Na = aN$,
and that $b \ne aN$ (otherwise $R(a,b)$ would be $I$), so
$N < aN < \fix(R(a,b)) < Q$, so $|Q:N| \ge 8$.
\end{proof}


\section{Subgroups of CC-loops}
\label{sec-sub}
Here, we show that some subloops of CC-loops are groups.

\begin{definition}
\label{def-assoc}
A triple of subsets $(A,B,C)$ of a loop $Q$ \emph{associates} iff
$x \cdot yz = xy \cdot z$ whenever
$x \in A$, $y \in B$, and $z \in C$.
A subset $S$ of $Q$ \emph{associates} iff
$(S,S,S)$ associates.
\end{definition}

Applying Theorem \ref{thm-perm},

\begin{lemma}
\label{lemma-assoc-perm}
In a CC-loop, the property ``$(A,B,C)$ associates'' is
invariant under all permutations of the set $\{A,B,C\}$.
\end{lemma}

By modifiying an argument of Bruck and Paige \cite{BP} for A-loops:

\begin{lemma}
In a CC-loop, if $(A,B,C)$ associates then
$(\sbl{A},\sbl{B},\sbl{C})$ associates.
\end{lemma}
\begin{proof}
For each $b \in B$ and $c \in C$, the map $R_b R_c R_{bc}\iv$ is
an automorphism (see Lemma \ref{lemma-aut})
and is the identity on $A$, so it is the identity
on $\sbl{A}$, which implies that
$(\sbl{A},B,C)$ associates.  
By Lemma \ref{lemma-assoc-perm}, we may apply this argument two
more times to prove that $(\sbl{A},\sbl{B},\sbl{C})$ associates.
\end{proof}

\begin{corollary}
\label{cor-assoc}
In a CC-loop, if $S$ associates, then
$\sbl{S}$ associates, and is hence a group.
\end{corollary}


\section{WIP Elements}
\label{sec-wip}

Throughout this section, $(Q, \cdot)$ always denotes a CC-loop.
By \cite{KUNC}, power-associative elements $x$ 
satisfy a number of additional properties. In this section,
we shall derive some further properties of these $x$
and their associated $E_x$ when $x$ is also a WIP element 
(see Definitions \ref{defn-E}, \ref{def-PA}, and \ref{def-WIP}).

Whenever $x$ is power-associative, 
all elements of the group generated by
$L_x$ and $R_x$ are of the form $E_x^r R_x^s L_x^t$ for
some $r,s,t\in\ZZZ$.  This is immediate from the following lemma,
which is taken from Lemmas 3.17 and 3.19 of \cite{KUNC}:

\begin{lemma}
\label{lemma-powers}
If $x$ is power-associative, then
for all $r,s,t,i,j,k,n\in \ZZZ$, the following hold:
\begin{itemizz}
\item[1.]
$E_x$ commutes with $L_x$ and $R_x$.
\item[2.]
$R_x^{-j} L_x^t R_x^j = E_x^{-jt} L_x^t$.
\item[3.]
$E_x^r R_x^s L_x^t \cdot E_x^i R_x^j L_x^k =
E_x^{r + i - jt} R_x^{s+j} L_x^{t+k}$.
\item[4.]
$R_{x^n} = E_x^{(n-1)n/2} R_x^n$
\item[5.]
$L_{x^n} = E_x^{-(n-1)n/2} L_x^n$.
\item[6.]
$E_{x^n} = E_x^{(n^2)}$.
\end{itemizz}
\end{lemma}

\begin{lemma}
\label{lemma-wip-powers}
In a CC-loop, if $c$ is a power-associative WIP element,
then for each $n\in\ZZZ$, $c^n$ is a WIP element.
\end{lemma}
\begin{proof}
Let $m = (n-1)n/2$.  Applying Lemma \ref{lemma-powers},
we have
\[
R_{c^n} \rho L_{c^n} = E_c^m R_c^n \rho E_c^{-m} L_c^n =
R_c^n \rho L_c^n = \rho \ \ .
\]
We are using the fact that $E_c$ commutes with $\rho$ (because
it is an automorphism) and with  $R_c$ (by Lemma \ref{lemma-powers}(1)).
\end{proof}

\begin{lemma}
\label{lemma-wip-eqns}
In a CC-loop, if $c$ is a power-associative WIP element,
then the following hold:
\[
\begin{array}{ccc}
D_c = L_{c\iv} \, \rho  &\qquad &  D_c\iv = R_{c\iv} \,  \lambda \\
\lambda L_c  \rho = R_c\iv &\qquad & \rho R_c \lambda = L_c\iv  \\
L_c \rho R_c = \rho  &\qquad & R_c \lambda L_c = \lambda 
\end{array}
\]
\end{lemma}

\begin{proof}
Note that since
$y D_c = y \ld c$ and $y D_c\iv = c \rd y$, the equations in the right
column are mirrors of the ones in the left, so we need only prove one
from each row.  For the first row, use
$ c \cdot g(c,y) z = y \cdot c z $
(see Lemma \ref{lemma-LR-conj}), and set $z = g(c,y)^\rho$ to get
$ c  = y \cdot c g(c,y)^\rho =  y \cdot c ( c\iv y \cdot c )^\rho $
(see Lemma \ref{lemma-getfg}).
Since $R_c \rho L_c = \rho$, we get
$ c  =   y \cdot (c\iv y)^\rho $, which implies $D_c = L_{c\iv}\,  \rho $.

For the second row, apply both equations in
the first row to $c\iv$, which is also WIP,
to get $L_c \rho = D_{c\iv} = \rho R_c \iv$.
The third row restates the second.
\end{proof}

Lemmas \ref{lemma-wip-equiv} and \ref{lemma-wip-eqns}
provide conjugation relations which, together with
Lemma \ref{lemma-powers}, show that
whenever $c$ is power-associative and WIP,
all elements of the group generated by
$L_c$, $R_c$, and $\rho$ are of the form $\alpha E_c^r R_c^s L_c^t$ for
some $r,s,t\in\ZZZ$, and some $\alpha \in \sbl{\rho}$.
It is also easy to see now that 
if $c$ is a power-associative WIP element,
then each of $R_c, L_c$ commutes with each of $\lambda^2, \rho^2$.

\begin{lemma}
\label{lemma-cc-mf}
In a CC-loop, $x (yz\cdot x) = (x^{\lambda}\ld y)\cdot zx$ and $(x\cdot yz) x = xy\cdot (z\rd x^{\rho})$.
\end{lemma}

\begin{proof}
By Lemmas \ref{lemma-CCAtop} and  \ref{lemma-getfg}
and Corollary \ref{coro-fginverses},
$(R_x, G_x, R_x)(F_x, L_x, L_x) = (L_{x^{\lambda}}\iv, R_x, R_x L_x)$
is an autotopism. Thus
$x (yz\cdot x) = (x^{\lambda}\ld y)\cdot zx$ for all $y,z$.
\end{proof}

\begin{lemma}
\label{lemma-wip-square}
In a CC-loop, if $c$ is a power-associative WIP element and 
$x$ is arbitrary, then
$x \cdot (x E_c\iv \cdot c) = x^2 \cdot c$.
\end{lemma}

\begin{proof}
We have $x^2 = (x^{\lambda}\ld c\iv) \cdot cx = (x\rd c) \cdot cx$
using Lemma \ref{lemma-cc-mf} and $D_{c\iv} = \rho R_c\iv$. Thus
$x^2 \cdot c = ((x\rd c) \cdot cx)\cdot c = (x\rd c)c\cdot g(c,cx)
= x\cdot xL_c L_{c\iv} R_c = x\cdot (x E_c\iv \cdot c)$ by
$\LCC$, Lemma \ref{lemma-getfg}, and Lemma \ref{lemma-powers}(5).
\end{proof}

\begin{lemma}
\label{lemma-wip-EE}
In a CC-loop, if $b$ and $c$ are power-associative, then
$\sbl{b,c}$ is a group if and only if $c E_b = c$ and $b E_c = b$.
If $c$ is also a WIP element, then $c E_b = c$ iff $b E_c = b$.
\end{lemma}

\begin{proof}
If $\sbl{b,c}$ is a group, then obviously $c E_b = c$ and $b E_c = b$.
Conversely, to prove that $\sbl{b,c}$ is a group, it is sufficient, by
Corollary \ref{cor-assoc}, to show that $\{b,c\}$ associates; that is,
$(x,y,z) = 1$ whenever $x,y,z \in \{b,c\}$. However, since $\sbl{c}$
and $\sbl{b}$ are groups and the associators are invariant under permutations
(Theorem \ref{thm-perm}), it is sufficient to show that $b^2 \cdot c = b \cdot bc$
and $c^2 \cdot b = c \cdot cb$. By Lemma \ref{lemma-powers}(5), these equations
are equivalent to $c E_b = c$ and $b E_c = b$, respectively.

Now if $c$ is a WIP element, then Lemmas \ref{lemma-wip-square}
and \ref{lemma-powers}(5) give
$b E_c\iv \cdot c = c L_{b^2} L_b\iv = b\cdot cE_b\iv$. Thus
$b E_c = b$ if and only if $c E_b = c$.
\end{proof}

\begin{lemma}
\label{lemma-wip-E2}
In a CC-loop, if $c$ is power-associative and WIP, then $E_c^2 = I$.
\end{lemma}

\begin{proof}
In Lemma \ref{lemma-cc-mf}, set $x = c$, $y = c\ld u$,
and $z = (c\ld u)^{\rho}$ to obtain
$c^2 = u ((c\ld u)^{\rho} / c\iv)$; equivalently, 
$D_{c^2} = L_c\iv \rho R_{c\iv}\iv$.
Now, applying Lemmas \ref{lemma-wip-eqns} and \ref{lemma-powers}, we get
$D_{c^2} = L_{c^{-2}} \, \rho = L_c^{-2} E_c^{-3} \rho$
and
$ L_c\iv \rho R_{c\iv}\iv =  L_c\iv  L_{c\iv} \rho  =  L_c^{-2} E_c^{-1} \rho$.
so that $ E_c^{-3} = E_c^{-1}$ .
\end{proof}

\begin{theorem}
\label{thm-group2}
In a CC-loop, if $c$ is WIP, and if $b$ and $c$ are power-associative, 
then $\sbl{b, c^2}$ is a group.
\end{theorem}

\begin{proof}
By Lemmas \ref{lemma-powers}(6) and \ref{lemma-wip-E2}, 
$b E_{c^2} = b E_c^4 = b$. Now apply Lemma \ref{lemma-wip-EE}
to $c^2$, which is WIP by Lemma \ref{lemma-wip-powers}.
\end{proof}

\begin{lemma}
\label{lemma-wip-Esquare}
In a CC-loop, if $c$ is WIP, and if $b$ and $c$ are power-associative, 
then $c E_b^2 = c$.
\end{lemma}

\begin{proof}
Lemma \ref{lemma-wip-square} implies $b\iv \cdot (b\ld (b^2\cdot c)) = 
b\iv \cdot (b E_c\iv \cdot c)$. Now $L_b^2 L_b\iv L_{b\iv} =
E_b^{-2}$ by Lemma \ref{lemma-powers}(5), and $b E_c\iv R_c
= b R_{c\iv}\iv = b\lambda D_c = b\iv \ld c$ by Lemmas
\ref{lemma-powers}(4) and \ref{lemma-wip-eqns}. Therefore
$c E_b^{-2} = b\iv \cdot (b\iv \ld c) = c$.
\end{proof}

\begin{theorem}
\label{thm-group1}
In a CC-loop, if $c$ is WIP, and if $b$ and $c$ are power-associative, 
then $\sbl{b^2, c}$ is a group.
\end{theorem}

\begin{proof}
By Lemmas \ref{lemma-powers}(6) and \ref{lemma-wip-Esquare}, 
$c E_{b^2} = c E_b^4 = c$. Now apply Lemma \ref{lemma-wip-EE}.
\end{proof}

Applying either Theorem \ref{thm-group2} or \ref{thm-group1}
we see that a power-associative WIP CC-loop in which every element is a square
must be a group.
Then, applying the Lagrange property (Corollary \ref{cor-lagrange}), we get:

\begin{corollary}
\label{cor-finite-square}
A finite power-associative WIP CC-loop of odd order is a group.
\end{corollary}

This corollary is not really new.
In \cite{BAS2}, {\cyr Basarab} shows that
a loop satisfies Wilson's
identity iff it is a ``generalized Moufang 
loop'' with squares in the nucleus. Then
Goodaire and Robinson \cite{GRB}
showed that a loop satisfies Wilson's
identity iff it is a WIP CC-loop.
Thus, in fact, all squares are nuclear in a WIP CC-loop,
so that $Q / N(Q)$ is a boolean group.
We give an example in Section \ref{sec-ex} of a power-associative WIP
CC-loop of order $16$ in which $|Q / N(Q)| = 4$; this is not 
an extra loop (that is, some $\sbl{b,c}$ fails to be a group),
so that Theorems \ref{thm-group2} and \ref{thm-group1} are best possible.


\section{Power-Associative CC-loops}
\label{sec-pow}

Throughout this section, $(Q, \cdot)$ always denotes a
power-associative CC-loop. We shall derive some further results
beyond Lemma \ref{lemma-powers}.
In particular, every cube is a WIP element (see Definition \ref{def-WIP}),
and each $E_x^6 = I$ (see Definition \ref{defn-E}).

\begin{lemma}
\label{lemma-LRD}
$R_x L_x = D_{x\iv} D_x$ and $L_x R_x = (D_x D_{x\iv})\iv$.
\end{lemma}

\begin{proof}
By Lemma \ref{lemma-getfg} and Corollary \ref{coro-fginverses},
$I = G_x F_x = D_{x\iv}\iv R_x \rho \cdot \rho  L_x D_x\iv =
D_{x\iv}\iv R_x  L_x D_x\iv$, and
$I = F_x G_x = D_{x\iv} L_x \lambda \cdot \lambda R_x D_x = 
D_{x\iv} L_x  R_x D_x$.
\end{proof}

Note that this lemma requires that power-associativity hold in $Q$,
not just that the particular element $x$ is power-associative,
since we needed $\rho^2 = I$, or equivalently, $\rho = \lambda$;
see Lemma \ref{lem-PA}.

Compare the following lemma with Corollary \ref{coro-preAAIP}.

\begin{lemma}
\label{lemma-short-B}
$(x \cdot xy)\cdot y\iv x\iv = x$ and
$x\iv y\iv\cdot  (yx \cdot x) = x$.
\end{lemma}

\begin{proof}
We compute
\[
\begin{array}{rcll}
L_x^2 &=&  E_x R_{x\iv} L_x^2 R_x &\qquad \mathrm{(Lemma\ \ref{lemma-powers})}\\
&=& E_x F_x D_{x\iv}\iv D_x\iv & \qquad \mathrm{(Lemmas\ \ref{lemma-getfg}\ and\ \ref{lemma-LRD})}\\
&=& E_x F_x G_x J R_x\iv D_x\iv & \qquad \mathrm{(Lemma\ \ref{lemma-getfg})} \\
&=& E_x J R_x\iv D_x\iv & \qquad \mathrm{(Corollary\ \ref{coro-fginverses})} \\
&=& J E_x R_x\iv D_x\iv & \qquad \mathrm{(Lemma\ \ref{lemma-aut})} \\
&=& J R_{x\iv} D_x\iv &\qquad \mathrm{(Lemma\ \ref{lemma-powers}(4))}\\
\end{array}
\]
Thus $x\cdot xy = x \rd (y\iv x\iv)$ or $(x\cdot xy) \cdot y\iv x\iv = x$, as claimed.
\end{proof}

\begin{lemma}
\label{lemma-short-C}
$y\iv \cdot (y R_x^3) = x^3$ and $(yL_x^3) \cdot y\iv = x^3$.
\end{lemma}

\begin{proof}
By Lemma \ref{lemma-LR-conj}, 
$ R_x\iv L_u R_x =  L_{x\iv  \ld u} L_x\iv $, so
$R_x L_{y\iv} = L_{x\iv y\iv} R_x L_x$.
Thus, $y\iv \cdot (y R_x^3) =
(yx\cdot x) R_x L_{y\iv} = (yx \cdot x) L_{x\iv y\iv}R_x L_x = x R_x L_x = x^3$
by Lemma \ref{lemma-short-B}.
\end{proof}

\begin{theorem}
\label{thm-cubes-WIP}
In a power-associative CC-loop, every cube is a WIP element.
\end{theorem}

\begin{proof} From Lemmas \ref{lemma-powers}(4,5) and \ref{lemma-short-C},
$E_x^3 L_{x^3} = L_x^3 = J D_{x^3}\iv$ and
$E_x^{-3} R_{x^3} = R_x^3 = J D_{x^3}$. Thus
$I = J E_x^3 L_{x^3} J E_x^{-3} R_{x^3} = J L_{x^3} J R_{x^3}$,
by Lemma \ref{lemma-powers}. Therefore $J L_{x^3} J = R_{x^3}\iv$,
that is, $x^3$ is a WIP element.
\end{proof}

\begin{corollary}
\label{cor-dias}
For each $b, c$ in a power-associative CC-loop, $\sbl{b, c^6}$ and
$\sbl{b^2, c^3}$ are groups.
\end{corollary}

\begin{proof}
$c^3$ is WIP, so apply Theorems \ref{thm-group2} and \ref{thm-group1}.
\end{proof}

The examples in Section \ref{sec-ex} show that some
$\sbl{b^2,c^2}$ can fail to be a group
(see Table \ref{table-27}), and so can some
$\sbl{b^3,c^3}$ (see Table \ref{table-16}).

\begin{corollary}
\label{cor-sixth}
If $Q$ is a power-associative CC-loop in which every element is
sixth power, then $Q$ is a group.
\end{corollary}

\begin{proof}
$Q$ is diassociative, and hence an extra loop.  However, in an extra loop, all
squares are in the nucleus \cite{FENB}, and so $Q = N(Q)$ is a group.
\end{proof}

Then, applying the Lagrange property (Corollary \ref{cor-lagrange}), we get:

\begin{corollary}
\label{cor-finite-sixth}
If $Q$ is a finite power-associative CC-loop of order
relatively prime to $6$, then $Q$ is a group.
\end{corollary}

In a power-associative CC-loop, Lemma \ref{lemma-powers}(6), 
Theorem \ref{thm-cubes-WIP}, and Lemma \ref{lemma-wip-E2}
imply $E_x^{18} = E_{x^3}^2 = I$.
We conclude this section with an improvement of this.

\begin{lemma}
\label{lemma-D2}
In a power-associative CC-loop, $x^2 = y \cdot ((x\iv y)\iv \cdot x)$ and $x^2 = (x\cdot (y x\iv)\iv)\cdot y$.
Thus $D_{x^2} = L_{x\iv} J R_x$ and $D_{x^2}\iv = R_{x\iv} J L_x$.
\end{lemma}

\begin{proof}
In Lemma \ref{lemma-cc-mf}, set $y = x\iv u$ and $z = (x\iv u)\iv$ to get
$x^2 = u \cdot ((x\iv u)\iv \cdot x)$.
\end{proof}

\begin{theorem}
\label{thm-E6}
Every power-associative CC-loop satisfies $E_x^6 = I$ for all $x$.
\end{theorem}

\begin{proof}
\[
\begin{array}{rcll}
E_x^{-3} L_x^6 &=& L_{x^2}^3 & \qquad \mathrm{(Lemma\ \ref{lemma-powers}(5))} \\
&=& J D_{x^6}\iv & \qquad \mathrm{(Lemma\ \ref{lemma-short-C})} \\
&=& J R_{x^{-3}} J L_{x^3} & \qquad \mathrm{(Lemma\ \ref{lemma-D2})} \\
&=& L_{x^{-3}}\iv L_{x^3} & \qquad \mathrm{(Theorem\ \ref{thm-cubes-WIP})} \\
&=& E_x^6 L_x^3 E_x^{-3} L_x^3 & \qquad \mathrm{(Lemma\ \ref{lemma-powers}(5))} \\
&=& E_x^3 L_x^6 & \qquad \mathrm{(Lemma\ \ref{lemma-powers}(1))}
\end{array}
\]
Rearranging, we have $E_x^6 = I$.
\end{proof}


\section{Semidirect Products}
\label{sec-semidir}

This standard construction from group theory generalizes 
to loops. We follow Goodaire and Robinson \cite{GRS}.

\begin{definition} \hskip 1cm 
\label{defn-semider}
\begin{itemizz}
\item[1.] Let $A, K$ be loops, and assume that $\varphi : A\to \Sym(K)$
satisfies $\varphi_{1_A} = I$ and $(1_K) \varphi_a = 1_K$ for
all $a\in A$. The \emph{external semidirect product} $A \ltimes_\varphi K$
is the set $A\times K$ with the binary operation
\[
(a, x)(b, y) := (a b,\, (x)\varphi_b  \cdot y) .
\]
for $a,b\in A$, $x,y\in K$. We write $A \ltimes K$ when $\varphi$
is clear from context.
\item[2.] A loop $Q$ is an \emph{internal semidirect product} of subloops
$A$ and $K$ if $K$ is normal in $Q$, $Q = AK$, $A\cap K = \{ 1 \}$,
and each of $(K,A,K)$, $(A,A,K)$, and $(A,K,Q)$ associates.
\end{itemizz}
\end{definition}

The external semidirect product $A\ltimes K$ is clearly a loop with left
and right division operations given, respectively, by
\[
\begin{array}{rcl}
(a,x)\ld (b,y) &=& (a \ld b,\, [(x)\varphi_{a\ld b}] \ld y )
\\
(a,x)\rd (b,y) &=& (a\rd b ,\, (x \rd y)\varphi_b\iv )
\end{array}
\]

The following comes from \cite{GRS}, Thms. 2.3 and 2.4.

\begin{proposition}\hskip 1cm  
\label{prop-semi}
\begin{itemizz}
\item[1.] If $Q = A\ltimes K$ is an external semidirect product of loops
$A$ and $K$, then $Q$ is isomorphic to the internal semidirect product
of the subloops $A\times \{ 1\}$ and $\{1\}\times K$.
\item[2.] If a loop $Q$ is an internal semidirect product of subloops
$A$ and $K$, then $Q$ is isomorphic to an external semidirect
product $A\ltimes_{\varphi} K$, where
$\varphi : A\to \Sym(K)$ is defined by:
$\varphi_a = R_a L_a\iv \rest K$.
\end{itemizz}
\end{proposition}

For CC-loops, the notion of semidirect product is much closer to its
group-theoretic specialization than for arbitrary loops. Recall from
Definition \ref{defn-nuclear} the notion of a nuclear automorphism.

\begin{lemma}
\label{lem-semiCC}
Let $Q$ be a CC-loop which is an internal semidirect product of
subloops $A$ and $K$, and define $\varphi : A\to \Sym(K)$ by
$\varphi_a := R_a L_a\iv \rest K$ for each $a\in A$. Then
$\varphi(A)\subseteq \NAut(K)$, and $\varphi : A\to \NAut(K)$
is a homomorphism.
\end{lemma}

\begin{proof}
Since $(A,K,Q)$ associates, we apply Lemma
\ref{lemma-assoc-perm} repeatedly in what follows
without explicit reference. (In CC-loops, the conditions
that $(K,A,K)$ and $(A,A,K)$ associate are redundant.)
For $x,y\in K$, $a\in A$,
\[
a\cdot (xy)\varphi_a = x\cdot ya = x(a \cdot (y)\varphi_a)
= xa\cdot (y)\varphi_a = (a\cdot (x)\varphi_a) \cdot (y)\varphi_a .
\]
Thus $(xy)\varphi_a = (x)\varphi_a \cdot (y)\varphi_a$, 
and so $\varphi_a\in \Aut(K)$. Now for each $x\in K$, $a\in A$,
Theorem \ref{thm-bas} implies there exists $c\in N(Q)$ such that
$x\varphi_a = xc$. But since $x\varphi_a \in K$, we have
$c\in K\cap N(Q)\subseteq N(K)$. Thus $\varphi_a\in \NAut(K)$.
Finally, for $a,b\in A$, $x\in K$,
we compute
\[
ab\cdot (x)\varphi_{ab} = xa\cdot b
= (a\cdot (x)\varphi_a)b = a ((x)\varphi_a\cdot b)
= a (b \cdot (x)\varphi_a\varphi_b). 
\]
Thus $(x)\varphi_{ab} = (x)\varphi_a\varphi_b$. This completes
the proof.
\end{proof}

We take notational advantage of Lemma \ref{lem-semiCC}
as follows: if $A\ltimes K$ is a CC-loop, then we set
$x^a := (x) \varphi_a$ for $x \in K$, $a \in A$. Note that
$x^{a^{\lambda}} = x^{a^{\rho}} = x\varphi_a\iv$.

We now prove that the necessary conditions for a semidirect
product to be a CC-loop given in Lemma \ref{lem-semiCC}
are also sufficient (Theorem \ref{thm-semider-CC}). This generalizes
D.A. Robinson's characterization of when $A\ltimes K$ is an extra
loop in the case where $A$ is a group \cite{ROB}.

\begin{theorem}
\label{thm-semider-CC}
Let $A, K$ be CC-loops, and $\varphi \in \Hom(A,\Aut(K))$.
Then the following are equivalent:
\begin{itemizz}
\item[1.] $A\ltimes_\varphi K$ is a CC-loop.
\item[2.] $\varphi_b\in \NAut(K)$ for all $b \in A$.
\item[3.] The triples
\[
\UU(x,b) := ( L_{x^b}R_x\iv, L_x , L_{x^b})
\qquad
and
\qquad
\VV(x,b) := (R_x,  R_{x^b} L_x\iv, R_{x^b})
\]
are in $\Atop(K)$ for all $x\in K$ and $b \in A$.
\end{itemizz}
\end{theorem}

\begin{proof}
For $(2) \leftrightarrow (3)$:
Fix $x\in K$ and $b \in A$. 
We have $\LL_x := (L_x R_x\iv , L_x, L_x) \in \Atop(K)$ by
Lemmas \ref{lemma-CCAtop} and \ref{lemma-getfg}.
Hence, $\UU(x,b) \LL_x\iv = ( L_{x^b} L_x\iv, I , L_{x^b} L_x\iv)$.
Now if $\UU(x,b)\in \Atop(K)$, then by
Lemma \ref{lemma-nuc-atop}(2), $1 L_{x^b} L_x\iv = x \ld (x^b) \in N(K)$,
so that $\varphi_b$ is a nuclear automorphism. Conversely, if $\varphi_b$ 
is nuclear, fix $k\in N(K)$ such that $x^b  = xk$.
Then $L_{x^b} L_x\iv = L_k$, and so 
$\UU(x,b) \LL_x\iv = ( L_k, I , L_k ) \in \Atop(K)$ by
Lemma \ref{lemma-nuc-atop}(1). Thus $\UU(x,b)\in \Atop(K)$ since $\Atop(K)$
is a group. A similar argument shows the equivalence of $\varphi_b\in \NAut(K)$
and $\VV(x,b)\in \Atop(K)$.

For $(1) \leftrightarrow (3)$:
Fix $(a,x), (b,y), (c,z)\in A\ltimes K$, and
write out the two sides of $\RCC$ in $A\ltimes K$ using
$f(u,v) = (uv)\rd u$ (Lemma \ref{lemma-getfg}) in $K$.
The left side is
\[
(a,x)\cdot (b,y)(c,z) =
(a\cdot bc,\; x^{bc} \cdot y^c  z)\ \  .
\]
The right side is
\[
[((a,x)(b,y))\rd (a,x)] \cdot (a,x)(c,z)
=
(f(a,b)\cdot ac,\;
[(x^{bc} y^c)\rd (x^c)]\cdot x^c z )\ \ .
\]
Equating the $K$-components, replacing $z$ by $z^c$
and then applying the automorphism $\varphi_c\iv$,
we get
$x^{b} \cdot y z = [(x^{b} y)\rd x] \cdot x z$.
Thus $A\ltimes K$ satisfies $\RCC$ iff each
$\UU(x,b)\in \Atop(K)$.

Likewise, we can write out the two sides of $\LCC$ in $A\ltimes K$ using
$g(u,v) = u \ld (vu)$ in $K$. The left side is
\[
(c,z)(b,y)\cdot (a,x) =
(cb\cdot a,\; z^{ba} y^a \cdot x)\ \  .
\]
The right side is
\[
(c,z)(a,x) \cdot [(a,x) \ld ((b,y)(a,x))]
=
(ca\cdot g(a,b),\;
(z^{ba} x^{a\ld (ba)}) [ x^{a\ld (ba)} \ld ( y^a  x )]  ) \ \ .
\]
Equating the $K$-components, replacing $x$ by $x^{b^{\lambda}a}$,
$z$ by $z^{b^{\lambda}}$, and then applying the automorphism
$\varphi_a\iv$, we get
$z y \cdot x^d =  z x \cdot [ x \ld ( y x^d )] $ where $d = b^{\lambda}$.
Thus $A\ltimes K$ satisfies $\LCC$ iff every $\VV(x,d)\in \Atop(K)$.
\end{proof}

We remark that the implication $(1)\rightarrow (2)$ follows directly
from Lemma \ref{lem-semiCC}. However, the proof of Theorem
\ref{thm-semider-CC} has the advantage of offering a
characterization of when $A\ltimes_\varphi Q$ satisfies $\RCC$ or
$\LCC$ alone, while our proof of Lemma \ref{lem-semiCC} relies on
Theorem \ref{thm-bas}. We also remark that in
proving $(1) \leftrightarrow (3)$, the arguments for $\LCC$ and $\RCC$
are similar, but we could not simply say that the $\LCC$ case follows from
the $\RCC$ case ``by mirror symmetry'', since there is an asymmetry in the
definition of $A\ltimes Q$.

Theorem \ref{thm-semider-CC} suggests that a natural definition of
\textit{holomorph} for a CC-loop $Q$ is $\NAut(Q) \ltimes_\varphi Q$, 
where $\varphi$ is the identity map. (This differs slightly from the usage
in \S5 of Bruck \cite{BR1}.)  If $Q$ is a group, then
$\NAut(Q) = \Aut(Q)$, and $\NAut(Q) \ltimes_\varphi Q$ reduces
to the usual definition of holomorph in group theory.


\section{Examples}
\label{sec-ex}

\begin{table}[htb]
{ 
\footnotesize
\arraycolsep=1.2pt
\newcommand\bl{\bullet}
\[
\begin{array}{c|ccccccccccccccccccccccccccc|}
\bl&  0& 1& 2& 3& 4& 5& 6& 7& 8& 9&10&11&12&13&14&15&16&17&18&19&20&21&22&23&24&25&26 \\
\hline
  0 &  0& 1& 2& 3& 4& 5& 6& 7& 8& 9&10&11&12&13&14&15&16&17&18&19&20&21&22&23&24&25&26 \\
  1 &  1& 2& 0& 4& 5& 3& 7& 8& 6&10&11& 9&13&14&12&16&17&15&19&20&18&22&23&21&25&26&24 \\
  2 &  2& 0& 1& 5& 3& 4& 8& 6& 7&11& 9&10&14&12&13&17&15&16&20&18&19&23&21&22&26&24&25 \\
  3 &  3& 4& 5& 6& 7& 8& 0& 1& 2&12&13&14&16&17&15&11& 9&10&22&23&21&24&25&26&20&18&19 \\
  4 &  4& 5& 3& 7& 8& 6& 1& 2& 0&13&14&12&17&15&16& 9&10&11&23&21&22&25&26&24&18&19&20 \\
  5 &  5& 3& 4& 8& 6& 7& 2& 0& 1&14&12&13&15&16&17&10&11& 9&21&22&23&26&24&25&19&20&18 \\
  6 &  6& 7& 8& 0& 1& 2& 3& 4& 5&15&16&17&11& 9&10&13&14&12&26&24&25&18&19&20&22&23&21 \\
  7 &  7& 8& 6& 1& 2& 0& 4& 5& 3&16&17&15& 9&10&11&14&12&13&24&25&26&19&20&18&23&21&22 \\
  8 &  8& 6& 7& 2& 0& 1& 5& 3& 4&17&15&16&10&11& 9&12&13&14&25&26&24&20&18&19&21&22&23 \\
  9 &  9&10&11&12&13&14&16&17&15&18&19&20&22&23&21&24&25&26& 0& 1& 2& 5& 3& 4& 8& 6& 7 \\
 10 & 10&11& 9&13&14&12&17&15&16&19&20&18&23&21&22&25&26&24& 1& 2& 0& 3& 4& 5& 6& 7& 8 \\
 11 & 11& 9&10&14&12&13&15&16&17&20&18&19&21&22&23&26&24&25& 2& 0& 1& 4& 5& 3& 7& 8& 6 \\
 12 & 12&13&14&15&16&17&10&11& 9&21&22&23&26&24&25&20&18&19& 4& 5& 3& 8& 6& 7& 1& 2& 0 \\
 13 & 13&14&12&16&17&15&11& 9&10&22&23&21&24&25&26&18&19&20& 5& 3& 4& 6& 7& 8& 2& 0& 1 \\
 14 & 14&12&13&17&15&16& 9&10&11&23&21&22&25&26&24&19&20&18& 3& 4& 5& 7& 8& 6& 0& 1& 2 \\
 15 & 15&16&17& 9&10&11&13&14&12&24&25&26&18&19&20&22&23&21& 8& 6& 7& 2& 0& 1& 3& 4& 5 \\
 16 & 16&17&15&10&11& 9&14&12&13&25&26&24&19&20&18&23&21&22& 6& 7& 8& 0& 1& 2& 4& 5& 3 \\
 17 & 17&15&16&11& 9&10&12&13&14&26&24&25&20&18&19&21&22&23& 7& 8& 6& 1& 2& 0& 5& 3& 4 \\
 18 & 18&19&20&21&22&23&26&24&25& 0& 1& 2& 5& 3& 4& 6& 7& 8& 9&10&11&13&14&12&16&17&15 \\
 19 & 19&20&18&22&23&21&24&25&26& 1& 2& 0& 3& 4& 5& 7& 8& 6&10&11& 9&14&12&13&17&15&16 \\
 20 & 20&18&19&23&21&22&25&26&24& 2& 0& 1& 4& 5& 3& 8& 6& 7&11& 9&10&12&13&14&15&16&17 \\
 21 & 21&22&23&24&25&26&20&18&19& 3& 4& 5& 6& 7& 8& 2& 0& 1&13&14&12&16&17&15& 9&10&11 \\
 22 & 22&23&21&25&26&24&18&19&20& 4& 5& 3& 7& 8& 6& 0& 1& 2&14&12&13&17&15&16&10&11& 9 \\
 23 & 23&21&22&26&24&25&19&20&18& 5& 3& 4& 8& 6& 7& 1& 2& 0&12&13&14&15&16&17&11& 9&10 \\
 24 & 24&25&26&18&19&20&23&21&22& 6& 7& 8& 1& 2& 0& 4& 5& 3&17&15&16&10&11& 9&14&12&13 \\
 25 & 25&26&24&19&20&18&21&22&23& 7& 8& 6& 2& 0& 1& 5& 3& 4&15&16&17&11& 9&10&12&13&14 \\
 26 & 26&24&25&20&18&19&22&23&21& 8& 6& 7& 0& 1& 2& 3& 4& 5&16&17&15& 9&10&11&13&14&12 \\
\hline
\end{array}
\]
}  
\caption{A Power-Associative CC-Loop}
\label{table-27}
\end{table}

\begin{table}[htb]
{ 
\footnotesize
\arraycolsep=1.2pt
\newcommand\bl{\bullet}
\[
\begin{array}{c|cccccccccccccccc|}
\bl&  0& 1& 2& 3& 4& 5& 6& 7& 8& 9&10&11&12&13&14&15 \\
\hline
  0 &  0& 1& 2& 3& 4& 5& 6& 7& 8& 9&10&11&12&13&14&15 \\
  1 &  1& 2& 3& 0& 5& 6& 7& 4& 9&10&11& 8&13&14&15&12 \\
  2 &  2& 3& 0& 1& 6& 7& 4& 5&10&11& 8& 9&14&15&12&13 \\
  3 &  3& 0& 1& 2& 7& 4& 5& 6&11& 8& 9&10&15&12&13&14 \\
  4 &  4& 5& 6& 7& 0& 1& 2& 3&12&13&14&15&10&11& 8& 9 \\
  5 &  5& 6& 7& 4& 1& 2& 3& 0&13&14&15&12&11& 8& 9&10 \\
  6 &  6& 7& 4& 5& 2& 3& 0& 1&14&15&12&13& 8& 9&10&11 \\
  7 &  7& 4& 5& 6& 3& 0& 1& 2&15&12&13&14& 9&10&11& 8 \\
  8 &  8& 9&10&11&15&12&13&14& 0& 1& 2& 3& 7& 4& 5& 6 \\
  9 &  9&10&11& 8&12&13&14&15& 1& 2& 3& 0& 4& 5& 6& 7 \\
 10 & 10&11& 8& 9&13&14&15&12& 2& 3& 0& 1& 5& 6& 7& 4 \\
 11 & 11& 8& 9&10&14&15&12&13& 3& 0& 1& 2& 6& 7& 4& 5 \\
 12 & 12&13&14&15&11& 8& 9&10& 6& 7& 4& 5& 3& 0& 1& 2 \\
 13 & 13&14&15&12& 8& 9&10&11& 7& 4& 5& 6& 0& 1& 2& 3 \\
 14 & 14&15&12&13& 9&10&11& 8& 4& 5& 6& 7& 1& 2& 3& 0 \\
 15 & 15&12&13&14&10&11& 8& 9& 5& 6& 7& 4& 2& 3& 0& 1 \\
\hline
\end{array}
\]
}  
\caption{A Power-Associative WIP CC-Loop}
\label{table-16}
\end{table}

The example in 
Table \ref{table-27} is a power-associative CC-loop of order $27$
and exponent three.
$Z(Q) = N(Q) = \{0,1,2\}$, and $\{0,1,2,3,4,5,6,7,8\}$
is a normal subloop.
Note that $|Z(Q)| = 3$ is required for non-associative CC-loops
of order $27$ by Corollary \ref{cor-center}.

This loop also has the Automorphic Inverse Property (AIP);
that is, $J\in\Aut(Q)$.

The example in Table \ref{table-16}
is a power-associative CC-loop of order $16$.
The loop must have the weak inverse property
because $3 \nmid 16$, so every element is a cube
(see Theorem \ref{thm-cubes-WIP}).
It is not diassociative because
$ 4 \cdot (8\cdot 4) \ne (4\cdot 8) \cdot 4$;
also, $Q = \sbl{4,8}$.
$\; Z(Q) = N(Q) = \{0,1,2,3\}$, so all squares are in the nucleus.

These examples were found by the program SEM \cite{ZZ}.
As usual, once one is given such an example,
it is easy to write a very short program (in, e.g., C or java or python)
to verify the claimed properties for it.



\begin{thebibliography}{99}

\bibitem{BAS2} {\cyr A. S. Basarab,
Ob odnom klasse} $G$-{\cyr lup},
{\itcyr Matematicheskie Issledovaniya }
{\cyr Tom 3, Vyp. 2 (8) }  (1968) 72--77.

\bibitem{BAS} {\cyr A. S. Basarab, Klass} LK-{\cyr lup},
{\itcyr Matematicheskie Issledovaniya } {\cyr Vyp. 120} (1991) 3--7.

\bibitem{BEL} {\cyr V. D. Belousov},
{\itcyr Osnovy Teorii Kvazigrupp i Lup},
{\cyr Izdatel\cprime stvo <Nauka>, Moskva}, 1967.

\bibitem{BR1} R. H. Bruck,
Contributions to the theory of loops,
\textit{Trans. Amer. Math. Soc.}  60 (1946) 245--354. 

\bibitem{BR2} R. H. Bruck,
{\it A Survey of Binary Systems}, Springer-Verlag, 1971.

\bibitem{BP} R. H. Bruck and L. J. Paige,
Loops whose inner mappings are automorphisms,
\textit{Ann. of Math.} (2)  (1956) 308--323.

\bibitem{CKRV} O.~Chein, M.~K.~Kinyon, A.~Rajah, and
P.~Vojt\v{e}chovsk\'y, Loops and the Lagrange property,
to appear in \textit{Results Math.}

\bibitem{FENA}  F. Fenyves, Extra loops I,
\textit{Publ. Math. Debrecen} 15 (1968) 235--238.

\bibitem{FENB}  F. Fenyves, Extra loops II,
\textit{Publ. Math. Debrecen} 16 (1969) 187--192.

\bibitem{GRA}  E. G. Goodaire and D. A. Robinson,
A class of loops which are isomorphic to all loop isotopes,
\textit{Canadian J. Math.} 34 (1982) 662--672.

\bibitem{GRB}  E. G. Goodaire and D. A. Robinson,
Some special conjugacy closed loops,
\textit{Canadian Math. Bull.} 33 (1990) 73--78.

\bibitem{GRS} E. G. Goodaire and D. A. Robinson,
Semi-direct products and Bol loops,
\textit{Demonstratio Math.} 27 (1994) 573--588.

\bibitem{KUNC} K. Kunen,
The structure of conjugacy closed loops,
\textit{Trans. Amer. Math. Soc.} 352 (2000) 2889--2911.

\bibitem{MC} W.W. McCune,
\textit{OTTER\ 3.0 Reference Manual and Guide},
Technical Report ANL-94/6, Argonne National Laboratory, 1994;
or see: \\
\verb+http://www-fp.mcs.anl.gov/division/software/+

\bibitem{OS} J. M. Osborn,
Loops with the weak inverse property,
\textit{Pacific J. Math.} 10 (1960) 295--304.

\bibitem{PF} H.~O.~Pflugfelder,
{\it Quasigroups and Loops: Introduction},
Sigma Series in Pure Math. \textbf{8}, Heldermann Verlag,
Berlin, 1990.

\bibitem{ROB} D.~A.~Robinson,
Holomorphy theory of extra loops,
\textit{Publ. Math. Debrecen} 18 (1971), 59--64.

\bibitem{SO} {\cyr L. R. So\u\i kis, O Spetsial\cprime nykh lupakh},
in {\itcyr Voprosy Teorii Kvazigrupp i Lup}
({\cyr V. D. Belousov }, ed.),
{\cyr Redakc.-Izdat. Otdel Akad. Nauk Moldav. SSR, Kishinev}, 1970,
pp.~122--131.

\bibitem{ZZ} J. Zhang and H. Zhang,
SEM: a system for enumerating models,
{\it Proc. 14th Int.. Joint Conf. on AI (IJCAI-95)},
Montr\'eal, 1995, pp. 298 -- 303; available at URL:
\verb+  http://www.cs.uiowa.edu/~hzhang/+

\end{thebibliography}
\end{document}